\documentclass{amsart}
\usepackage{amsmath,amssymb,amscd,amsthm,hyperref}
\usepackage[all]{xy}
\makeatletter
\newtheorem*{rep@theorem}{\rep@title}
\newcommand{\newreptheorem}[2]{
	\newenvironment{#1}[1]{\def\rep@title{#2 \ref{##1}}\begin{rep@theorem}}{\end{rep@theorem}}}
\makeatother

\makeatletter
\def\thm@space@setup{
  \thm@preskip=\parskip \thm@postskip=0pt
}
\makeatother

\newcommand{\fm}{\mathfrak{m}}
\newcommand{\fp}{\mathfrak{p}}

\newcommand{\bF}{\mathbb{F}}
\newcommand{\bZ}{\mathbf{Z}}

\newcommand{\cD}{\mathcal{D}}

\newcommand{\cL}{\mathcal{L}}
\newcommand{\cM}{\mathcal{M}}
\newcommand{\cN}{\mathcal{N}}
\newcommand{\cO}{\mathcal{O}}

\newcommand{\cQ}{\mathcal{Q}}

\DeclareMathOperator{\ev}{ev}

\DeclareMathOperator{\Hom}{Hom}
\DeclareMathOperator{\End}{End}
\DeclareMathOperator{\Ext}{Ext}

\DeclareMathOperator{\Spec}{Spec}

\DeclareMathOperator{\colim}{colim}

\DeclareMathOperator{\Sh}{Sh}

\DeclareMathOperator{\Crys}{Crys}
\DeclareMathOperator{\Loc}{Loc}

\newcommand{\et}{\textrm{\'{e}t}}

\newcounter{ctr} \numberwithin{ctr}{section}

\theoremstyle{definition}
\newtheorem{dfn}[ctr]{Definition}

\theoremstyle{remark}
\newtheorem{rmk}[ctr]{Remark}

\theoremstyle{plain}
\newtheorem{thm}[ctr]{Theorem}
\newtheorem*{thm*}{Theorem}
\newtheorem{prp}[ctr]{Proposition}
\newtheorem*{prp*}{Proposition}
\newtheorem{lem}[ctr]{Lemma}

\newreptheorem{rthm}{Theorem}
\newreptheorem{rprp}{Proposition}

\newcommand{\stacks}[1]{{\href{http://stacks.math.columbia.edu/tag/#1}{[Stacks #1]}}}

\begin{document}

\title{Tannakian properties of unit Frobenius-modules}
\author{Maxim Mornev}
\email{m.mornev@math.leidenuniv.nl}
\address{Mathematisch Instituut \\
Universiteit Leiden \\
Postbus 9512 \\
2300 RA Leiden \\
Nederland 
}\address{
Dipartimento di Matematica ``Federigo Enriques'' \\
Universit\`{a} degli Studi di Milano \\
via Cesare Saldini 50 \\
20133 Milano MI \\
Italia
}
\address{Korteweg--de Vries Instituut voor Wiskunde\\
Universiteit van Amsterdam \\
Postbus 94248 \\
1090 GE Amsterdam \\
Nederland}
\begin{abstract}
We show that unit $\cO_{F,X}^\Lambda$-modules of Emerton and Kisin
provide an analogue of locally constant sheaves in the context of
B\"{o}ckle--Pink $\Lambda$-crystals. For example they form a tannakian
category if the coefficient algebra $\Lambda$ is a field. Our results
hold for a big class of coefficien algebras which includes Drinfeld
rings, and for arbitary locally noetherian base schemes.
\end{abstract}
\maketitle

\section*{Introduction}

Let $q$ be a prime power, $\bF_q$ a field with $q$ elements. For an
$\bF_q$-scheme $X$, and a commutative $\bF_q$-algebra $\Lambda$,
B\"{o}ckle and Pink \cite{bp} introduced the category of
$\Lambda$-crystals on $X$ which is in many ways analogous to the
category of \'{e}tale constructible sheaves of $\Lambda$-modules. Let us
recall the definitions.

\begin{dfn}[B\"{o}ckle--Pink \cite{bp}]\label{basic} Let $X$ be an
$\bF_q$-scheme, $\Lambda$ a commutative $\bF_q$-algebra. Let $F$ be the
endomorphism of $\Spec \Lambda \times_{\bF_q} X$ which acts as identity
on $\Spec \Lambda$, and as the absolute $q$-Frobenius on $X$. Let
$\cO_X^\Lambda$ denote the structure sheaf of $\Spec \Lambda
\times_{\bF_q} X$. 

(1) An \emph{$\cO_{F,X}^\Lambda$-module}
is a pair $(\cM,\varphi)$, where $\cM$ is an $\cO_X^\Lambda$-module,
and $\varphi\colon F^\ast \cM \to \cM$ a morphism of
$\cO_X^\Lambda$-modules. A morphism of $\cO_{F,X}^\Lambda$-modules
$\alpha\colon(\cM,\varphi_\cM) \to (\cN,\varphi_\cN)$ is a morphism
$\alpha\colon \cM \to \cN$ of $\cO_X^\Lambda$-modules, such that $\alpha
\circ \varphi_\cM = \varphi_\cN \circ F^\ast \alpha$.
The category of $\cO_{F,X}^\Lambda$-modules which are quasi-coherent as
$\cO_X^\Lambda$-modules is denoted $\mu(X,\Lambda)$.

(2) An $\cO_{F,X}^\Lambda$-module $(\cM,\varphi)$ is called
\emph{nilpotent} if for $n \gg 0$ the composition $\varphi \circ F^\ast
(\varphi) \circ \dotsc \circ F^{\ast n}(\varphi)$ vanishes. A morphism of
$\cO_{F,X}^\Lambda$-modules is called a \emph{nil-isomorphism} if its
kernel, and cokernel are nilpotent.

(3) The category $\Crys(X,\Lambda)$ of \emph{$\Lambda$-crystals} on $X$
is the localization of the category of $\cO_{F,X}^\Lambda$-modules which
are $\cO_X^\Lambda$-coherent at the multiplicative system of
nil-isomorphisms.\end{dfn}

The connection with constructible sheaves is provided by the
following result of B\"{o}ckle and Pink. Let $\Sh(X_\et,\Lambda)$
be the category of \'{e}tale sheaves of $\Lambda$-modules. Define a
functor $\varepsilon\colon \mu(X,\Lambda) \to \Sh(X_\et,\Lambda)$ by
\begin{equation*}
\varepsilon(\cM,\varphi)(u\colon U \to X) =
\Hom_{\mu(U,\Lambda)}\Big((\cO_U^\Lambda,1), (u^\ast\cM,u^\ast\varphi)\Big).
\end{equation*}
Since $\varepsilon$ transforms nil-isomorphisms to
isomorphisms (\cite{bp} propostion 10.1.7 (b)) one gets a
functor $\varepsilon\colon \Crys(X,\Lambda) \to \Sh(X_\et,\Lambda)$. 
\begin{thm}[B\"{o}ckle--Pink \cite{bp} Theorem 10.3.6]\label{bpeq}
Assume that $\Lambda$ is of finite dimension as $\bF_q$-vector space, 
$X$ is noetherian, and separated over $\bF_q$. The functor $\varepsilon$
defines an equivalence of $\Crys(X,\Lambda)$, and the category of
constructible \'{e}tale sheaves of $\Lambda$-modules.
\end{thm}
In particular $\varepsilon$ indentifies the subcategory $\Loc(X,\Lambda)
\subset \Sh(X_\et,\Lambda)$
of locally constant sheaves with a certain subcategory of
$\Crys(X,\Lambda)$.  An explicit description of this subcategory is
provided by a classical result of Katz which in the language of this
text can be stated as follows.
\begin{dfn}[Emerton--Kisin \cite{ek} Definition 5.1] An
$\cO_{F,X}^\Lambda$-module $(\cM,\varphi)$ is called \emph{unit} if
$\varphi$ is an isomorphism. The category of unit
$\cO_{F,X}^\Lambda$-modules which are locally of finite type as
$\cO_X^\Lambda$-modules is denoted $U(X,\Lambda)$.\end{dfn}
\begin{thm}[Katz \cite{katz} Proposition 4.1.1]\label{katzeq}
Assumptions as in theorem \ref{bpeq}.

\emph{(1)} The natural functor $U(X,\Lambda) \to \Crys(X,\Lambda)$ is
fully faithful.

\emph{(2)} The functor $\varepsilon$ identifies the subcategories
$U(X,\Lambda)$ and $\Loc(X,\Lambda)$.\end{thm}

The definition of $\Lambda$-crystals makes sense even if $\Lambda$ is not finite
over $\bF_q$. The interest in infinite coefficient algebras $\Lambda$
comes from a construction of Drinfeld (see e.g. section 3.5 of
\cite{bp}) which produces an $\cO_{F,X}^\Lambda$-module, and a fortiori
a $\Lambda$-crystal out of a sufficiently good Drinfeld module
$\varphi\colon \Lambda \to \End_{\bF_q}(\cL)$, $\cL$ a line bundle on $X$.
So one has a natural source of $\Lambda$-crystals which have
arithmetic significance.

For infinite $\Lambda$ the connection with \'{e}tale sheaves breaks down:
the functor $\varepsilon$ may assing a zero sheaf to a nonzero
crystal. Nevertheless B\"{o}ckle and Pink show that under mild technical
assumptions on $X$, and $\Lambda$, the crystals retain many properties
one would expect from constructible sheaves. Crystals form an abelian
category. Some of the Grothendieck six operations are defined for them:
pullback, proper pushforward, tensor product, and extension by zero.
These operations are related in the same way as they are for
constructible sheaves, and the equivalence of theorem \ref{bpeq} is
compatible with them. The pullback is an exact operation, and a crystal
is zero if and only if its pullbacks to all points $x \to
X$ are. 

Thus one can view $\Lambda$-crystals as an extension, or a modification
of the notion of a constructible sheaf. From such a point of view it is
natural to ask which kind of crystals is the analogue of locally
constant sheaves? In this paper we purport to show that the answer is
unit modules.

\begin{rmk} The part (1) of theorem \ref{katzeq} holds for infinite
$\Lambda$ as well provided one restricts to unit modules which are
$\cO_X^\Lambda$-coherent. So unit modules are indeed a particular kind of
crystals. Since we will not use this property we only indicate a
proof. A cohomological computation in the spirit of \cite{la} section 4
shows that if $\cM$ is a unit module, and $\cN$ is nilpotent then
$\Ext^1(\cM,\cN) = 0$ in $\mu(X,\Lambda)$.  Since $\Hom(\cM,\cN)$ also
vanishes it follows that a nil-isomorphism to a unit module is
necessarily a split epimorphism whence the result.\end{rmk}

If $\Lambda$ is finite-dimensional over $\bF_q$, $X$ is connected, and
locally noetherian then a choice of a geometric point $\overline{x} \to
X$ identifies $\Loc(X,\Lambda)$ with the category of
$\pi_1^{\textrm{\'{e}t}}(X,\overline{x})$-representations in
$\Lambda$-modules. In particular if $\Lambda$ is a finite field
extension of $\bF_q$ then $\mathrm{Loc}(X,\Lambda)$ is a tannakian
$\Lambda$-linear category.  It is this last property which we take as
the characteristic property of local systems. For infinite $\Lambda$ the
connection of $U(X,\Lambda)$ with the \'{e}tale fundamental group is
lost. Nevertheless in the spirit of B\"{o}ckle and Pink \cite{bp} we
prove the following:

\begin{rthm}{tan}Suppose that $\Lambda$ is a field, $X$ is connected, and
locally noetherian. The category $U(X,\Lambda)$ is
$\Lambda$-linear tannakian with respect to a rigid monoidal
structure inherited from $\Crys(X,\Lambda)$.\end{rthm}

It seems that no good notion of a tannakian
$\Lambda$-linear category is known in the case when $\Lambda$ is not a
field. Still in this case we demonstrate that under mild technical
assumptions on $\Lambda$, and $X$, unit modules retain some important tannakian
properties which local systems have. Unit modules form an abelian
category (theorem \ref{cohu}), and the pullback functors on them are
faithful
(theorem \ref{faithpull}). A nice feature of these results is that they
hold not only for regular but for arbitrary locally noetherian schemes
$X$.

In order to deduce these theorems we study structural properties of
locally finitely generated $\cO_X^\Lambda$-modules which admit a unit module
structure. Our approach was inspired by the following folklore
lemma
: if a $\cD$-module is $\cO_X$-coherent then it is
locally free.  It has a counterpart in the setting of Emerton--Kisin:
\begin{prp}[Emerton--Kisin \cite{ek} proposition 6.9.3]\label{ekrig}
Assume that $\Lambda$ is finite-dimensional as $\bF_q$-vector space, $X$
is locally noetherian. If a unit $\cO_{F,X}^\Lambda$-module is
$\cO_X$-coherent then it is locally free.
\end{prp}
We extend this proposition to the case of infinite $\Lambda$ in two
ways. One of them works for arbitrary $X$, and $\Lambda$.
\begin{rthm}{katzflat} If a unit $\cO_{F,X}^\Lambda$-module is locally
of finite presentation as an $\cO_X^\Lambda$-module then it is
$\cO_X$-flat.\end{rthm}
The main difficulty in the proof of this theorem is that when $\Lambda$
is infinite unit modules are almost never of finite type as
$\cO_X$-modules. Another extension of proposition \ref{ekrig} is
specific for the case when $\Lambda$ is a field:

\begin{rprp}{ufinpres} Assume that $\Lambda$ is a field, $X$ is locally
noetherian. Let $\cM$ be a unit $\cO_{F,X}^\Lambda$-module. If $\cM$ is
locally of finite type as an $\cO_X^\Lambda$-module then it is a unit
$(\Lambda,F)$-crystal (\cite{ek}, definition 6.9.1), i.e. it is locally
free as an $\cO_X^\Lambda$-module. \end{rprp}

This proposition follows easily from the next result which we call the
invariant closed subscheme theorem:

\begin{rthm}{invsubsch} Assume that $X$ is connected, and locally noetherian. Let $Z
\subset \Spec \Lambda \times_{\bF_q} X$ be a closed subscheme. If
$F^{-1} Z = Z$ as a subscheme then $Z = Z_0 \times_{\bF_q} X$ where
$Z_0$ is the scheme-theoretic image of $Z$ in $\Spec \Lambda$.\end{rthm}
Its proof is based on the proof of lemma 4.6.1 from \cite{bp}.

We leave several important questions unanswered.  We do not know how to
describe the Tannaka-dual groups (or gerbes) of the tannakian categories
which we construct. How do they change when the coefficient algebra
varies? If one fixes the base scheme $X$, and varies $\Lambda$ then our
unit modules form a stack. It is tempting to call it a stack of
tannakian categories. Is there a kind of Tannaka-dual object 
behind the whole stack?

\begin{rmk} The notion of unit modules comes from the work of Emerton and
Kisin \cite{ek}. However we use unit modules in a way which is different
from \cite{ek}. Emerton and Kisin construct a characteristic $p > 0$
version of the Riemann--Hilbert correspondence. In their work,
$\cO_{F,X}^\Lambda$-modules are analogues of $\cD$-modules, and unit
modules play the role of holonomic $\cD$-modules. In this text we work
not on the $\cD$-module side of Emerton--Kisin but on the side of
constructible sheaves, i.e. $\Lambda$-crystals.\end{rmk}
 
\subsection*{Notation and conventions}

Throughout the text $\Lambda$ denotes a commutative $\bF_q$-algebra
which we will use as the coefficients for unit modules.  For an
$\bF_q$-scheme $X$ the symbol $F$ indicates the endomorphism of $\Spec
\Lambda \times_{\bF_q} X$ which acts as identity on $\Spec \Lambda$, and
as the absolute $q$-Frobenius on $X$. We make no general assumptions on
the coefficient algebra $\Lambda$, and the base scheme $X$. In each
theorem we carefully state the precise assumptions on $\Lambda$, and $X$
under which it holds.

Following \cite{ek} we denote $\mu(X,\Lambda)$ the category of
$\cO_{F,X}^\Lambda$-modules whose underlying $\cO_X^\Lambda$-module is
quasi-coherent. 

$U(X,\Lambda)$ stands for the subcategory of $\mu(X,\Lambda)$ consisting
of unit modules whose underlying $\cO_X^\Lambda$-modules are locally of
finite type. We do not use the adjective ``coherent'' since in the
situation of interest for us it may well happen that $\Spec \Lambda
\times_{\bF_q} X$ is not a locally noetherian scheme.

\subsection*{Acknowledgements} Our debt to
Gebhard B\"{o}ckle,
Matthew Emerton,
Nicholas M. Katz,
Mark Kisin,
Vincent Lafforgue,
and Richard Pink
is obvious from the text. Besides them we would like to thank Lenny
Taelman who attracted our attention to this circle of problems, and
whose constant support and encouragement made the present work possible.
We would also like to thank Fabrizio Andreatta for several fruitful
discussions. The Stacks Project \cite{s} was of immense help in the
course of this work.
Our research is supported by a Leiden/Huygens fellowship, and an
ALGANT doctoral scholarship.

\section{The invariant closed subscheme theorem}\label{sinvsubsch}

\begin{thm}\label{invsubsch} Let $X$ be a connected locally noetherian
scheme over $\bF_q$, $\Lambda$ an $\bF_q$-algebra, and $Z \subset \Spec
\Lambda \times_{\bF_q} X$ a closed subscheme. If $F^{-1} Z = Z$ then $Z
= Z_0 \times_{\bF_q} X$ where $Z_0$ is the scheme-theoretic image of $Z$
under the projection $\Spec \Lambda \times_{\bF_q} X \to \Spec
\Lambda$.\end{thm}

The proof of this theorem is based on the proof of lemma 4.6.1 from \cite{bp}.

\begin{rmk}\label{perfectoid} When $X$ is not locally noetherian theorem
\ref{invsubsch} fails even if $\Lambda = \bF_q$. For an example let $R$
be the $q$-perfection of $\bF_q[[t]]$. The
ring $R$ is a local domain with maximal ideal $\fm$ generated by
$\{t^{q^n}\}_{n \in \bZ}$. Since $F$ maps this system of generators to
itself the ideal $\fm$ is $F$-invariant, and defines a closed
$F$-invariant subscheme of $\Spec R$ which is different from $\Spec R$,
and $\varnothing$.\end{rmk}

{\it Proof of theorem \ref{invsubsch}.}
(1) Assume that $X$ is the
spectrum of a local noetherian ring $R$ with maximal ideal $\fm$. Let $I
\subset \Lambda \otimes_{\bF_q} R$ be the ideal of $Z$.  The ideal $I$
is $F$-invariant in the sense that it is generated by $F(I)$. Consider
the ideal $I \cap \Lambda$ where $\Lambda$ is viewed as a subring of
$\Lambda \otimes_{\bF_q} R$ via the coprojection. Our goal is to show
that the inclusion $(I \cap \Lambda) \otimes_{\bF_q} R \subset I$ is an
equality.

It is enough to consider the case $I \cap \Lambda = 0$. Indeed, let
$\Lambda_0 = \Lambda/(I \cap \Lambda)$, and let $J$ be the image of $I$ in $\Lambda_0
\otimes_{\bF_q} R$. Then $J$ is $F$-invariant, and $J \cap \Lambda_0 =
0$, so the theorem implies that $J = 0$, i.e. $I = (I \cap \Lambda)
\otimes_{\bF_q} R$.

In order to deduce the equality $I = 0$ from $I \cap \Lambda = 0$ we only
need to prove the inclusion $I \subset \Lambda \otimes_{\bF_q} \fm$
since $F$-invariance of $I$ will then imply the inclusion $I \subset
\Lambda \otimes_{\bF_q} \bigcap_{n \geqslant 0} \fm^{q^n} = 0$.
Choose a basis $\{\lambda_s\}_{s \in S}$ of $\Lambda$ as an
$\bF_q$-vector space. It induces a basis of $\Lambda \otimes_{\bF_q} R$
as an $R$-module. In this basis $F$ acts by raising coordinates to
$q$-th powers.

Suppose that there are elements $f \in I$ which have a coordinate
belonging to $R^\times$. Among them pick one which has minimal number of
nonzero coordinates, say $f = r_1 \lambda_1 + r_2 \lambda_2 + \dotsc r_n
\lambda_n$ with $r_1 \in R^\times$, and all $r_i \ne 0$. Dividing $r_1$
out we may assume that the coordinate of $f$ at $\lambda_1$ is $1$.
The element $F(f) - f \in I$ has coefficient
zero at $\lambda_1$, so by minimality of $n$ the elements $r_i^q - r_i$
belong to $\fm$. Since $r_i^q - r_i = \prod_{\alpha \in \bF_q} (r_i -
\alpha)$, and since $\fm$ is prime we conclude that for each $i$ there
exist $\alpha_i \in \bF_q$, $m_i \in \fm$ such that $r_i = \alpha_i +
m_i$.

Next consider the set $I_1 \subset I$ of elements having the form
$\lambda_1 + (\alpha_2 + m_2)\lambda_2 + \dotsc + (\alpha_n +
m_n)\lambda_n$ for some $\alpha_i \in \bF_q, m_i \in \fm$.  Pick an
element $f \in I_1$, and let $i \in \{1 \dotsc n\}$ be the maximal index
with the property that $m_j = 0$ for all $j \leqslant i$. Consider the element
$F(f) - f \in I$. Its coordinates at
$\lambda_j, j \leqslant i,$ are zero, and the coordinate at
$\lambda_{i+1}$ is $m_{i+1}^q - m_{i+1}$ which is equal to $m_{i+1}$ up
to a unit. Dividing by this unit, and subtracting the result from $f$ we
obtain an element of $I_1$, which has $m_j = 0$ for all $j \leqslant i +
1$. Repeating this process we see that $I_1$ contains an element of
$\Lambda$. As all elements of $I_1$ are nonzero this contradicts the
assumption $\Lambda \cap I = 0$. Hence $I \subset \Lambda
\otimes_{\bF_q} \fm$ as we want.

(2) Let $X$ be a connected locally noetherian scheme. Pick a point $x
\in X$. According to the step (1) the base change of $Z$ to $\Spec \Lambda
\times_{\bF_q} \Spec \cO_{X,x}$ has the form $Z_x \times_{\bF_q} \Spec
\cO_{X,x}$ for a certain closed subscheme $Z_x \subset \Spec \Lambda$.

Let $x, x' \in X$ be points such that $x$ is a specialization of $x'$.
There exists a discrete valuation ring $R$ and a morphism $f\colon \Spec
R \to X$ mapping the generic point of $R$ to $x'$, and the
closed point to $x$ \stacks{054F}. Step (1) shows that the base change
of $Z$ to $\Spec \Lambda \times_{\bF_q} \Spec R$ along $1_\Lambda
\times_{\bF_q} f$ has the form $Z_f \times_{\bF_q} \Spec R$ for a
certain $Z_f \subset \Spec \Lambda$. Therefore $Z_x = Z_f = Z_{x'}$.
Since every two points of $X$ can be connected by a chain of
generalizations, and specializations, we see that $Z_x = Z_0$ is
independent of $x \in X$.
As the ideal sheaves of $Z$, and $Z_0 \times_{\bF_q} X$ have the same
stalks at every point of $\Spec\Lambda \times_{\bF_q} X$ these closed
subschemes are equal. \qed

We will use theorem \ref{invsubsch} through the following corollary:

\begin{lem}\label{coeffred} Let $X$ be a connected locally noetherian
scheme over $\bF_q$, $\Lambda$ an $\bF_q$-algebra, $\cM$ a unit
$\cO_{F,X}^\Lambda$-module which is locally of finite type as an
$\cO_X^\Lambda$-module. Let $Z_n(\cM) \subset \Spec \Lambda
\times_{\bF_q} X$ be the closed subscheme defined by the $n$-th Fitting
ideal sheaf of $\cM$ as an $\cO_X^\Lambda$-module. The subscheme
$Z_n(\cM)$ has the form $Z \times_{\bF_q} X$ for a certain closed
subscheme $Z \subset \Spec \Lambda$.\end{lem}

{\it Proof.} $F$-invariance of $\cM$ implies $F$-invariance of $Z_n(\cM)$,
whence the result.\qed

\begin{rmk} Here are two more easy but interesting corollaries of
theorem \ref{invsubsch}. Let $X$ be a connected locally noetherian
scheme over $\bF_q$, $\Lambda$ an $\bF_q$-algebra, and $\cM$ a unit
$\cO_{F,X}^\Lambda$-module which is locally of finite type as an
$\cO_X^\Lambda$-module.

(1) If for every maximal ideal $\fm \subset \Lambda$
the restriction of $\cM$ to the closed subscheme $\Spec \Lambda/\fm
\times_{\bF_q} X \subset \Spec \Lambda \times_{\bF_q} X$ is zero then
$\cM = 0$.

(2) If a fiber of $\cM$ over a point $x \in X$ is locally free of rank
$r$ as an $\cO^\Lambda_{k(x)}$-module then $\cM$ is locally free of rank
$r$ as an $\cO^\Lambda_X$-module.

What makes (1) a nontrivial statement is the fact that in general there
are closed subsets of $\Spec \Lambda \times_{\bF_q} X$ whose image in
$\Spec \Lambda$ misses all the closed points.  \end{rmk}

\section{Flatness}\label{sflat}

In this section we will prove that under suitable finiteness assumptions
every unit $\cO_{F,X}^\Lambda$-module is $\cO_X$-flat (theorem
\ref{katzflat}). We will deduce this theorem from the following 
lemma.

\begin{lem}\label{flartin} Let $R$ be an artinian local ring over
$\bF_q$, $\fm$ its maximal ideal, $k$ the residue field.

\emph{(1)} Equip $k$ with the structure of an $\bF_q$-algebra via the
composite arrow $\bF_q \to R \to k$. Every splitting $k \to R$ of the
quotient map $R \to k$ provided by Cohen structure theorem is a morphism
of $\bF_q$-algebras.

\emph{(2)} Fix a splitting as in (1). Let $\Lambda$ be an
$\bF_q$-algebra. If $M$ is a unit $\cO_{F,R}^\Lambda$-module which is of
finite presentation as an $\cO_R^\Lambda$-module then $M \cong
(M/\fm) \otimes_k R$ over $\cO_R^\Lambda$.\end{lem}

\begin{rmk} Let $U_{pf}(X,\Lambda) \subset U(X,\Lambda)$ be the
subcategory of modules which are of finite presentation as
$\cO_X^\Lambda$-modules. In fact under the assumptions of lemma
\ref{flartin} (2) the reduction functor
\begin{equation*}
-\otimes_R R/\fm\colon U_{pf}(R,\Lambda)
\to U_{pf}(k,\Lambda)
\end{equation*}is an equivalence of categories, its inverse being
$-\otimes_k R$. One can say that local artinian rings have a henselian
property with respect to unit modules. We do not know if other types of
rings exhibit such a property.\end{rmk}

{\it Proof of lemma \ref{flartin}.} (1) Let $x \in \bF_q \subset R$, and
pick $y \in k \subset R$ which reduces to the same element as $x$ in
$k$. Set $z = x - y$.  Then $z^q = z$, so $(z^{q-1})^2 = z^q \cdot
z^{q-2} = z^{q-1}$. Since $R$ is local it follows that $z^{q-1} \in \{ 0, 1 \}$.
By assumption $z$ is not a unit, so $z = z \cdot z^{q-1} = 0$. 

(2) Let $A$ be a matrix which defines a presentation of $M$ as an
$\Lambda \otimes_{\bF_q} R$-module, $F^n(A)$ the matrix obtained from
$A$ by applying $F$ componentwise $n$ times. $F^n(A)$ defines a
presentation of $F^{\ast n} M$, and hence of $M$. Since $R$ is local
artinian the coefficients of $F^n(A)$ lie in the subring $\Lambda
\otimes_{\bF_q} k \subset \Lambda \otimes_{\bF_q} R$ for $n$ big enough,
whence the claim. \qed

\begin{thm}\label{katzflat} Let $X$ be a scheme over $\bF_q$, $\Lambda$
an $\bF_q$-algebra. If $(\cM,\varphi)$ is a unit
$\cO_{F,X}^\Lambda$-module which is locally of finite presentation as an
$\cO_X^\Lambda$-module then $\cM$ is flat over $\cO_X$.\end{thm}

\begin{rmk}
We expect that a stronger result is true. Namely, if $X =
\Spec R$ is a local noetherian $\bF_q$-algebra, $\Lambda$ an arbitrary
$\bF_q$-algebra, and $M$ a unit $\cO_{F,X}^\Lambda$-module of finite
presentation as an $\cO_X^\Lambda$-module then $M$ is free as an
$R$-module.

Even in the case $\Lambda = \bF_q$ theorem \ref{katzflat}
fails if the module $\cM$ is not of finite presentation over
$\cO_X^\Lambda$. For an example let $R$ be the ring of remark
\ref{perfectoid}, $\fm$ its maximal ideal. $F$-invariance of $\fm$ implies
that $R/\fm$ has a structure of a unit
$\cO_{F,R}^{\bF_q}$-module, and this module is not $R$-flat.
\end{rmk}

{\it Proof of theorem \ref{katzflat}.} The question being local on $X$
it is enough to consider the case $X = \Spec R$. Let $M$ be the
$\Lambda\otimes_{\bF_q} R$-module of corresponding to $\cM$.

(1) Assume that $R$ is noetherian, and $\Lambda$ is of finite type over
$\bF_q$. Let $N_1 \to N_2$ be an inclusion of $R$-modules of finite
type, $K = \ker(M \otimes_R (N_1 \to N_2))$. We want to show that $K =
0$. As $K$ is of finite type over $S = \Lambda \otimes_{\bF_q} R$ it is
enough to show that $K \otimes_R k(\fp) = 0$ for every prime $\fp \in
\Spec R$.

So we may assume $R$ is local. Let $\fm$ be its maximal ideal,
$\widehat{S}$ the completion of $S$ at $\fm$, and $\widehat{M} = M
\otimes_S \widehat{S}$. For every ideal $I \subset R$ the module $I
\otimes_R \widehat{M}$ is $\fm$-adically separated as it is the
completion of $I \otimes_R M$. Moreover for all $n > 0$ the
$R/\fm^n$-module $\widehat{M}/\fm^n = M/\fm^n$ is free by lemma
\ref{flartin}. Hence $\widehat{M}$ is $R$-flat by
theorem 1 of \cite[ch. III, \S5]{bo}. So
$K \otimes_S \widehat{S} = 0$, and $K/\fm = 0$.

(2) Let $\Lambda$, $R$ be arbitrary. By \stacks{05N7} there exist
finitely generated $\bF_q$-subalgebras $\Lambda_0 \subset
\Lambda$, $R_0 \subset R$, a finitely presented $\Lambda_0
\otimes_{\bF_q} R_0$-module $M_0$, and morphisms $f\colon F^\ast M_0 \to
M_0$, $g\colon M_0 \to F^\ast M_0$ such that $M \cong \Lambda
\otimes_{\Lambda_0} \otimes M_0 \otimes_{R_0} R$, and under this
isomorphism the structure morphism $\varphi\colon F^\ast M \to M$ is
the base extension of $f$, $\varphi^{-1}$ is the base extension of
$g$, and moreover $f g = 1$, $gf = 1$. Thus $(M_0, f)$ is a 
unit $\cO_{F,R_0}^{\Lambda_0}$-module.

Represent $\Lambda$ as a filtered colimit of finitely generated
$\Lambda_0$-subalgebras $\Lambda_i$. If $I \subset R$ is an ideal then
\begin{equation*}
\Lambda \otimes_{\Lambda_0} M_0 \otimes_{R_0} (I \to R) = 
\colim_i \big(\Lambda_i \otimes_{\Lambda_0} M_0 \otimes_{R_0} (I \to R)
\big)
\end{equation*}
The module $\Lambda_i \otimes_{\Lambda_0} M_0$ is $R_0$-flat by (1), so
the morphism on the left hand side being a filtered colimit of
injections is an injection itself. \qed

\section{Tannakian properties}\label{stan}

\begin{thm}\label{cohu} Let $X$ be a locally noetherian scheme over
$\bF_q$, $\Lambda$ an algebra essentially of finite type over $\bF_q$.
The category $U(X,\Lambda)$ is closed under kernels, and cokernels in
$\mu(X,\Lambda)$. In particular, it is an abelian category.\end{thm}

The only obstruction this theorem has to deal with is the 
non-exactness of $F^\ast$. According to Kunz \cite{kunz} $F^\ast$ is
exact if and only if $X$ is regular. So the main interest in theorem
\ref{cohu} comes from the fact that it holds for singular base schemes
$X$ too.

{\it Proof of theorem \ref{cohu}.} Let $\alpha\colon (\cM,\varphi_\cM) \to
(\cN,\varphi_\cN)$ be a morphism of unit modules which is an epimorphism of
underlying $\cO_X^\Lambda$-modules, and let $(\cQ,\varphi_\cQ)$ be the
kernel of $\alpha$ computed in $\mu(X,\Lambda)$. The lower row of the
diagram
\begin{equation*}
\xymatrix{
0 \ar[r] &
\cQ \ar[r] &
\cM \ar[r]^{\alpha} &
\cN \ar[r] &
0 \\
0 \ar[r] &
F^\ast\cQ \ar[r] \ar[u]_{\varphi_\cQ} &
F^\ast\cM \ar[r]^{F^\ast \alpha} \ar[u]_{\varphi_\cM} &
F^\ast\cN \ar[r] \ar[u]_{\varphi_\cN} &
0
}
\end{equation*}
is exact because $\cN$ is $\cO_X$-flat by theorem \ref{katzflat}.
Thus $\varphi_\cQ$ is an isomorphism. Using epi-mono factorization in
$\mu(X,\Lambda)$, and the fact that $U(X,\Lambda)$ is closed under
cokernels one concludes that $\ker \alpha \in U(X,\Lambda)$ for general
$\alpha$.\qed

If $f\colon X \to Y$ is a morphism of schemes then pullback of the
underlying $\cO_X^\Lambda$-module defines a functor $f^\ast\colon
U(X,\Lambda) \to U(Y,\Lambda)$.

\begin{thm}\label{faithpull} Let $f\colon Y \to X$ be a morphism of
locally noetherian schemes over $\bF_q$, $\Lambda$ an algebra
essentially of finite type over $\bF_q$. 

\emph{(1)} The pullback functor $f^\ast\colon
U(X,\Lambda) \to U(Y,\Lambda)$
is exact.

\emph{(2)} If $X$ is connected, and $Y \ne \varnothing$ then $f^\ast$ is
faithful, and conservative.\end{thm}

{\it Proof.} (1) Indeed according to theorem \ref{katzflat} the 
modules in question are $\cO_X$-flat.

(2) If $\cM \in U(X,\Lambda)$ then by lemma
\ref{coeffred} $Z_0(\cM) = Z \times_{\bF_q} X$ for a certain closed
subscheme $Z \subset \Spec \Lambda$. Hence $Z_0(f^\ast \cM) = Z
\times_{\bF_q} Y$ can be empty if and only if $Z = \varnothing$, i.e.
$\cM = 0$. \qed

\begin{dfn}\label{tensor} Let $X$ be an $\bF_q$-scheme, $\Lambda$ an
$\bF_q$-algebra.

(1) The tensor product of two unit modules $(\cM,\varphi_\cM)$, $(\cN,\varphi_\cN) \in
U(X,\Lambda)$ is defined as
\begin{equation*}
(\cM,\varphi_\cM) \otimes (\cN,\varphi_\cN) = (\cM \otimes_{\cO_X^\Lambda}
\cN, \varphi_\cM \otimes_{\cO_X^\Lambda} \varphi_\cN),
\end{equation*}
and similarly for morphisms.

(2) The symbol $1$ denotes the object $(\cO_X^\Lambda, 1) \in
U(X,\Lambda)$.

(3) The constraints of associativity, commutativity, and unity for this
tensor product are taken from the respective constraints for
$\cO_X^\Lambda$-modules.
\end{dfn}
It is easy to check that the tensor product above provides
$U(X,\Lambda)$ with a symmetric monoidal structure.

\begin{dfn}\label{dual} Let $X$ be an $\bF_q$-scheme, $\Lambda$ an
$\bF_q$-algebra, and $(\cM,\varphi) \in U(X,\Lambda)$. Assume that
$\cM$ is locally free as an $\cO_X^\Lambda$-module.

(1) Let $(\cM,\varphi)^\vee$ be the module $(\cM^\vee,
(\varphi^\vee)^{-1})$ where dualization happens with respect to the
$\cO_X^\Lambda$-module structure.

(2) The evaluation $\ev\colon (\cM,\varphi) \otimes (\cM,\varphi)^\vee
\to 1$, and coevaluation $\delta\colon 1 \to (\cM,\varphi)^\vee \otimes
(\cM,\varphi)$ morphisms are lifted from the corresponding morphisms of
locally free $\cO_X^\Lambda$-modules.
\end{dfn}
It is straightforward to check that the evaluation, and
coevaluation morphisms of (2) are morphisms of
$\cO_{F,X}^\Lambda$-modules. They satisfy identities (2.1.2) of
\cite{d} by construction.
%

In the case when $\Lambda$ is a field unit modules enjoy an important
property.
 
\begin{prp}\label{ufinpres} Let $X$ be a locally noetherian scheme over
$\bF_q$, $\Lambda$ a field containing $\bF_q$. Every $\cM \in
U(X,\Lambda)$ is a unit $(\Lambda,F)$-crystal (\cite{ek}, definition
6.9.1), i.e. it is locally free as an $\cO_X^\Lambda$-module.\end{prp}

{\it Proof.} We can assume $X$ is connected. Since $\Lambda$
is a field lemma \ref{coeffred} shows that Fitting subschemes $Z_n(\cM)$
are either empty or of the form $\Spec \Lambda \times_{\bF_q} X$. Hence
$\cM$ is locally free \stacks{07ZD}.\qed

Thus the dual objects of definition \ref{dual} exist for every unit
module, and $U(X,\Lambda)$ becomes a rigid monoidal category. Also note
that under assumptions of the proposition $\cO_X^\Lambda$ is coherent in
the sense of \stacks{01BV} so every $\cM \in U(X,\Lambda)$ is a
coherent $\cO_X^\Lambda$-module.

\begin{thm}\label{tan} Let $X$ be a connected locally noetherian scheme
over $\bF_q$, $\Lambda$ a field containing $\bF_q$. The category
$U(X,\Lambda)$ equipped with the tensor structure above is
$\Lambda$-linear tannakian, and the inclusion $U(X,\Lambda) \subset
\mu(X,\Lambda)$ is exact. If $X(\bF_q) \ne \varnothing$ then
$U(X,\Lambda)$ is neutral, i.e.  has a fiber functor to $\Lambda$-vector
spaces.\end{thm}

A nice feature of this theorem is that it places no restriction on the
size of $\Lambda$ as an $\bF_q$-algebra. One can take $\Lambda =
\bF_q((t))$ or the completion of its algebraic closure.
Note that abelianness of $U(X,\Lambda)$ does not
follow directly from theorem \ref{cohu}.

{\it Proof of theorem \ref{tan}.} We will verify that
$(U(X,\Lambda),\otimes)$ is a tensor category as
defined in \cite{d}, section 2. It is then tannakian since pullback
to a point $y \in \Spec \Lambda \times_{\bF_q} X$ provides a fiber
functor to $k(y)$-vector spaces. If $s\colon
\Spec \bF_q \to X$ is a section of the structure morphism $X \to \Spec
\bF_q$ then the pullback functor $s^\ast\colon
U(X,\Lambda) \to
U(\bF_q,\Lambda)$ gives rise to a fiber functor with
values in $\Lambda$-vector spaces.

The conditions (2.1.1), (2.1.2) of \cite{d} mean precisely that the
category in question is rigid symmetric monoidal which we know from
proposition \ref{ufinpres}.

(2.1.3). $U(X,\Lambda)$ is closed under cokernels in
$\mu(X,\Lambda)$. Therefore in order to check that it is
closed with respect to kernels it is enough to consider a morphism $\alpha\colon
(\cM,\varphi_\cM)~\to~(\cN,\varphi_\cN)$ which is an epimorphism of underlying
$\cO_X^\Lambda$-modules. Let $(\cQ,\varphi_\cQ)$ be the kernel of $\alpha$
computed in $\mu(X,\Lambda)$.
The modules $\cM$, $\cN$ are locally free by proposition \ref{ufinpres}.
Hence they are locally of finite presentation over $\cO_X^\Lambda$, and
as a consequence $\cQ$ is locally of finite type \stacks{0519}. Since
$\cN$ is $\cO_X^\Lambda$-flat it follows that $\varphi_\cQ$ is an
isomorphism.

(2.1.4). We need to check that the natural morphism $\Lambda \to
\End(1)$ is bijective. The endomorphism ring in question consists of $f
\in \Gamma(\Spec\Lambda \times_{\bF_q} X, \cO_X^\Lambda)$ which are
invariant under $F$. Since $\Gamma(\Spec\Lambda \times_{\bF_q} X,
\cO_X^\Lambda) = \Lambda \otimes_{\bF_q} \Gamma(X, \cO_X)$ the claim
follows from the theory of Artin-Schreier equation.\qed

\end{document}